\documentclass[a4paper,12pt]{article}
\usepackage{amsthm,amssymb,amsmath}

\title{Some Thoughts on Multiparameter Stochastic Processes}
\author{Niels Jacob and Alexander Potrykus}

\begin{document}
\maketitle
In this paper we want to set the stage for a new approach to multiparameter processes. The principle idea is to find a ``substitute'' for the notion ``generator'' as known from the theory of Markov processes. Thus our investigations are more analytic in their nature and will deal with the problem of constructing multiparameter processes starting with analytic data. Eventually we are lead to rather ``strange'' classes of (pseudo-)differential equations. These equations in general do not fit to any traditional class -- however they enter naturally in our context.\\
\\
The excellent book \cite{KHO} of D. Khoshnevisan and the references given therein may serve to get an overview on the existing theory. This theory handles multiparameter processes under the assumption that marginal processes are Markov and that the multiparameter process is given by ``commuting'' marginal processes (semigroups, generators).\\
\\
Let us start with some simple observations. For simplicity we will work in the beginning with only two time parameters $s,t\in\mathbb{R}_+$.\\
\\
Consider a continuous function $b:\mathbb{R}_+\times \mathbb{R}_+\times \mathbb{R}^n\rightarrow\mathbb{C}$ such that $b(s,t;0)=0$ for all $s,t\geq 0$ and $\xi\mapsto b(s,t;\xi)$ is a continuous negative definite function, compare \cite{BER} or \cite{JAC1} for a definition and characterization. It follows by the easy part of the Schoenberg theorem, compare Chr. Berg and G. Forst \cite{BER}, p. 41, that for every $s,t\geq 0$ there exists a probability measure $\mu_{(s,t)}$ on the Borel sets of $\mathbb{R}^n$ characterized by
\begin{align}\label{eq1}
\widehat{\mu}_{(s,t)}(\xi):=(2\pi)^{-\frac{n}{2}}e^{-b(s,t;\xi)}\quad,\xi\in\mathbb{R}^n.
\end{align}
Thus we can consider $\mu_{(s,t)}$ as distribution of a random variable $X_{(s,t)}:\Omega\rightarrow \mathbb{R}^n$. Clearly we can construct a large product space on which all the random variables $X_{(s,t)}$, $s,t\geq 0$, are defined. This however need not lead to an interesting process.\\
\\
Consider the following example. Let $\psi:\mathbb{R}^n\rightarrow\mathbb{C}$, $\psi(0)=0$ be a fixed continuous negative definite function and put
\begin{align}\label{eq2}
\widehat{\mu}_{(s,t)}(\xi)=(2\pi)^{-\frac{n}{2}}e^{-ts^2\psi(\xi)}.
\end{align}
If we restrict $(s,t)\in\mathbb{R}_+^2$ to a ``curve'' $(s,t)\mapsto (\sqrt{s},t_0)$, $t_0>0$ fixed, we find
\begin{align}\label{eq3}
\widehat{\mu}_{(\sqrt{s},t_0)}(\xi)=(2\pi)^{-\frac{n}{2}}e^{-t_0 s \psi(\xi)},
\end{align}
and if we restrict $(s,t)\in\mathbb{R}_+^2$ to a ``curve'' $(s,t)\mapsto (s_0,t)$, $s_0>0$ fixed, we get
\begin{align}\label{eq4}
\widehat{\mu}_{(s_0,t)}(\xi)=(2\pi)^{-\frac{n}{2}}e^{-s_0 t\psi(\xi)}.
\end{align}
Thus in both cases we have associated L\'evy processes. More precisely, the family of probability measures $(\mu_{(s,t)})_{s,t\geq 0}$ gives rise to two families of L\'evy processes. Thus, although there is no obvious way how to associate a process with $(\mu_{(s,t)})_{s,t\geq 0}$, often we can find two families of processes which should serve as substitute for the marginal processes of the classical approach.\\
\\
Let us remark that in the above considerations nothing will change if we substitute $b:\mathbb{R}^2_+\times\mathbb{R}^n\rightarrow \mathbb{C}$ by a function
$b:\mathbb{R}^k_+\times\mathbb{R}^n\rightarrow\mathbb{C}$ which we require to
be continuous and for which $\xi\mapsto b(s_1,\ldots s_k;\xi)$ is a continuous
negative definite function with $b(s_1,\ldots,s_k;0)=0$, i.e. we consider the measures given by
\begin{align}\label{eq5}
\widehat{\mu}_{(s_1,\ldots,s_k)}(\xi)=(2\pi)^{-\frac{n}{2}}e^{-b(s_1,\ldots,s_k;\xi)}.
\end{align}
Now, in some cases we may encounter the situation that there are $k$ functions $g_j$ such that for every $1\leq j\leq k$
\begin{align}\label{eq5a}
\widehat{\mu}_{(s_1^0,\ldots,g_j(s_j),\ldots,s_k^0)}(\xi)=(2\pi)^{-\frac{n}{2}}e^{-b(s_1^0,\ldots,g_j(s_j),\ldots,s_k;\xi)}
\end{align}
gives rise to a convolution semigroup, hence a L\'evy process.
On the Schwartz space $\mathcal{S}(\mathbb{R}^n)$ we may define the operators
\begin{align}\label{eq6}
T_{(s_1,\ldots,s_k)}u(x):=(2\pi)^{-\frac{n}{2}}\int_{\mathbb{R}^n} e^{i(x,\xi)}
e^{-b(s_1,\ldots,s_k,\xi)}\widehat{u}(\xi)\,\mathrm{d} \xi.
\end{align}
Since $\xi\mapsto e^{-b(s_1,\ldots,s_k;\xi)}$ is positive definite,
$$
\big|e^{-b(s_1,\ldots,s_k;\xi)}\big|\leq e^{-b(s_1,\ldots,s_k;0)}=1
$$
by our assumption that $b(s_1,\ldots,s_k;0)=0$, and therefore \eqref{eq6} is well defined. Clearly we have on $C_{\infty}(\mathbb{R}^n)$ with the obvious extension to $B_b(\mathbb{R}^n)$
\begin{align}\label{eq7}
T_{(s_1,\ldots,s_k)}u(x)=\int_{\mathbb{R}^n} u(x-y)\mu_{(s_1,\ldots,s_k)}(\,\mathrm{d} y)
\end{align}
with
\begin{align}\label{eq8}
\widehat{\mu}_{(s_1,\ldots,s_k)}(\xi)=(2\pi)^{-\frac{n}{2}}e^{-b(s_1,\ldots,s_k;\xi)}.
\end{align}
It follows that each of the operators $T_{(s_1,\ldots,s_k)}$ are contractions on $L^2(\mathbb{R}^n)$ as well as on $C_{\infty}(\mathbb{R}^n)$ and $B_b(\mathbb{R}^n)$, respectively. As translation invariant operators they also commute, i.e. for $s_1,\ldots,s_k\in\mathbb{R}_+$ and $t_1,\ldots,t_k\in \mathbb{R}_+$ it holds
\begin{align}\label{eq9}
\big[T_{(s_1,\ldots,s_k)},T_{(t_1,\ldots,t_k)}\big]=0.
\end{align}
Now, let us add the assumption that $(s_1,\ldots,s_k)\mapsto b(s_1,\ldots,s_k;\xi)$ is $k$-times continuously differentiable and that for $\sigma\in\mathbb{N}_0^k$, $|\sigma|\leq k$, it holds
\begin{align}\label{eq10}
\big|\partial_{s_1,\ldots,s_k}^{\sigma}b(s_1,\ldots s_k;\xi)\big|\leq C_{\sigma}\big(1+|\xi|^2\big)^{\frac{r}{2}}
\end{align}
where $r=r(\sigma)$. In this case we are allowed to perform the following calculation for each $u\in\mathcal{S}(\mathbb{R}^n)$
\begin{align*}
\frac{\partial^k}{\partial s_1\ldots\partial s_k}T_{(s_1,\ldots,s_k)}u(x)&=
\frac{\partial^k}{\partial s_1\ldots\partial s_k}(2\pi)^{-\frac{n}{2}}\int_{\mathbb{R}^n} e^{i(x,\xi)}e^{-b(s_1,\ldots, s_k;\xi)}\widehat{u}(\xi)\,\mathrm{d} \xi\\
&=(2\pi)^{-\frac{n}{2}}\int_{\mathbb{R}^n} e^{i(x,\xi)}\Bigg(\frac{\partial^k}{\partial s_1\ldots\partial s_k}e^{-b(s_1,\ldots,s_k;\xi)}\Bigg)\widehat{u}(\xi)\,\mathrm{d}\xi\\
&=(2\pi)^{-\frac{n}{2}}\int_{\mathbb{R}^n} e^{i(x,\xi)}a(s_1,\ldots,s_k;\xi)e^{-b(s_1,\ldots,s_k;\xi)}\widehat{u}(\xi)\,\mathrm{d} \xi
\end{align*}
with
\begin{align}\label{eq11}
a(s_1,\ldots,s_k;\xi)=e^{b(s_1,\ldots,s_k;\xi)}\frac{\partial^k}{\partial s_1\ldots\partial s_k}e^{-b(s_1,\ldots,s_k;\xi)}.
\end{align}
Denoting by $a(s_1,\ldots,s_k,D_x)$ the pseudodifferential operator 
\begin{align}\label{eq12}
a(s_1,\ldots,s_k,D_x)u(x)=(2\pi)^{-\frac{n}{2}}\int_{\mathbb{R}^n} e^{i(x,\xi)}a(s_1,\ldots,s_k;\xi)\widehat{u}(\xi)\,\mathrm{d} \xi
\end{align}
we find that for $u\in\mathcal{S}(\mathbb{R}^n)$ a solution to the equation
\begin{align*}
\frac{\partial^k}{\partial s_1\ldots\partial s_k}v(s_1,\ldots,s_k;x)=a(s_1,\ldots,s_k;D_x)v(s_1,\ldots,s_k;x)
\end{align*}
is given by
\begin{align*}
v(s_1,\ldots,s_k;x):=T_{(s_1,\ldots,s_k)}u(x).
\end{align*}
Moreover, for $1\leq j\leq k$ it follows that
\begin{align}\label{eq13}
\lim_{s_j\rightarrow 0}v(s_1,\ldots,s_k;x)=(2\pi)^{-\frac{n}{2}}\int_{\mathbb{R}^n} e^{i(x,\xi)}e^{-b(s_1,\ldots,s_{j-1},0,s_{j+1},\ldots,s_k;\xi)}\widehat{u}(\xi)\,\mathrm{d}\xi.
\end{align}
For $k=2$ let us investigate the scope of operators $a(s_1,s_2,D_x)$ for which
we will write now $a(s,t;D)$ only.
\begin{enumerate}
\item Consider with two continuous negative definite functions $\psi_j:\mathbb{R}^n\rightarrow\mathbb{C}$ the function
\begin{align}\label{eq14}
(s,t)\mapsto s\psi_1(\xi)+t\psi_2(\xi).
\end{align}
This gives
\begin{align}\label{eq15}
a(s,t;\xi)=\psi_1(\xi)\psi_2(\xi).
\end{align}
In particular, for $\psi_1(\xi)=\psi_2(\xi)=|\xi|^2$ we find
\begin{align}\label{eq16}
a(s,t;D_x)=\Delta_x^2.
\end{align}
Moreover, with $n=1$, $\psi_1(\xi)=\xi^2$ and $\psi_2(\xi)=ic\xi$, $c\in\mathbb{R}$, we get $\psi_1(\xi)\psi_2(\xi)=ic\xi^3$ which corresponds to the operator
\begin{align}\label{eq17}
a(s,t;D_x)=-c\frac{\partial^3}{\partial x^3}.
\end{align}
Finally, with $n=1$, $\psi_1(\xi)=i\xi$ and $\psi_2(\xi)=-i\xi$ we have $\psi_1(\xi)\psi_2(\xi)=\xi^2$ and obtain the operator
\begin{align}\label{eq18}
a(s,t;D_x)=-\Delta_x.
\end{align}
\item We decompose $\mathbb{R}^n=\mathbb{R}^{n_1}\times\mathbb{R}^{n_2}$ and consider on $\mathbb{R}^{n_j}$ the continuous negative definite function
$\psi_j:\mathbb{R}^{n_j}\rightarrow\mathbb{C}$. Clearly the function
$$
(s,t;\xi)\mapsto b(s,t;\xi)=s\psi_1(\xi^1)+t\psi(\xi^2)\:,\xi^j\in\mathbb{R}^{n_j}, \xi=(\xi^1,\xi^2),
$$
is for $s,t\geq 0$ fixed a continuous negative definite function and for the
corresponding symbol $a(s,t;\xi)$ we find
\begin{align}\label{eq19}
a(s,t;\xi)=\psi_1(\xi^1)\psi_2(\xi^2)
\end{align}
which gives
\begin{align}\label{eq20}
a(s,t;D)=\psi_1(D_{x_1})\circ \psi_2(D_{x_2})\quad,x^j\in\mathbb{R}^{n_j},x=(x^1,x^2)
\end{align}
For $\psi_j(\xi^j)=|\xi^j|^2$ we arrive at
\begin{align}\label{eq21}
a(s,t;D)=\Delta_{n_1}\Delta_{n_2}
\end{align}
This operator was investigated by E. B. Dynkin \cite{DYN}, note also the preceeding work of J. B. Walsh (and co-authors), see \cite{WAL} and the references therein.
\item Let us now consider the case where with two continuous negative definite
functions $\psi_j:\mathbb{R}^n\rightarrow\mathbb{C}$ and a continuous function $d:\mathbb{R}_+^2\times\mathbb{R}^n\rightarrow \mathbb{C}$, such that
$d(0,t;\xi)=d(s,0;\xi)=d(s,t;0)=0$ and $\xi\mapsto d(s,t;\xi)$ is negative
definite, the function $b$ is given by
\begin{align}\label{eq22}
b(s,t;\xi)=s\psi_1(\xi)+t\psi_2(\xi)+d(s,t;\xi).
\end{align}
A straightforward calculation leads to
\begin{align}\label{eq23}
\begin{split}
a(s,t;\xi)=&\psi_1(\xi)\psi_2(\xi)+\psi_1(\xi)\frac{\partial}{\partial t}d(s,t;\xi)+\psi_2(\xi)\frac{\partial}{\partial s}d(s,t;\xi)\\
&\mspace{50mu}+\frac{\partial}{\partial s}d(s,t;\xi)\frac{\partial}{\partial t}d(s,t;\xi)-\frac{\partial^2}{\partial s\partial t}d(s,t;\xi).
\end{split}
\end{align}
From our general considerations we find for $\varphi\in\mathcal{S}(\mathbb{R}^n)$ that
\begin{align}\label{eq24}
\mspace{-10mu}v(s,t;\xi)\mspace{-3mu}:=\mspace{-3mu}T_{(s,t)}\varphi(x)\mspace{-5mu}=\mspace{-5mu}(2\pi)^{-\frac{n}{2}}\int_{\mathbb{R}^n}\mspace{-15mu} e^{i(x,\xi)}e^{-s\psi_1(\xi)-t\psi_2(\xi)-d(s,t;\xi)}\widehat{\varphi}(\xi)\,\mathrm{d}\xi
\end{align}
solves
\begin{align}\label{eq25}
\frac{\partial^2}{\partial s\partial t} v(s,t;x)=a(s,t;D_x)v(s,t;x).
\end{align}
In this case we can also find the boundary or initial conditions
\begin{align}\label{eq26}
\lim_{s\rightarrow 0}T_{(s,t)}\varphi(x)=(2\pi)^{-\frac{n}{2}}\int_{\mathbb{R}^n} e^{i(x,\xi)}e^{-t\psi_2(\xi)}\widehat{\varphi}(\xi)\,\mathrm{d}\xi=T_t^{(2)}\varphi(x)
\intertext{and} \label{eq27}
\lim_{t\rightarrow 0}T_{(s,t)}\varphi(x)=(2\pi)^{-\frac{n}{2}}\int_{\mathbb{R}^n} e^{i(x,\xi)}e^{-s\psi_1(\xi)}\widehat{\varphi}(\xi)\,\mathrm{d} \xi=T_s^{(1)}\varphi(x),
\end{align}
where $\big(T_s^{(1)}\big)_{s\geq 0}$ is the Feller semigroup associated with $\psi_1$ and $\big(T_t^{(2)}\big)_{t\geq 0}$ is the Feller semigroup corresponding to $\psi_2$.
Thus by \eqref{eq24} a solution to the following problem is given:
\begin{quote}
Given $\varphi\in\mathcal{S}(\mathbb{R}^n)$ and extend $\varphi$ by $\big(T_s^{(1)}\big)_{s\geq 0}$ as well as by $\big(T_t^{(2)}\big)_{t\geq 0}$ to $\mathbb{R}_+\times \{0\}\times\mathbb{R}^n$ and $\{0\}\times\mathbb{R}_+\times\mathbb{R}^n$, respectively. Then $T_{(s,t)}$, $s,t\geq 0$, is the solution operator, better the two-parameter family of solution operators, to \eqref{eq25}--\eqref{eq27}.
\end{quote}
In some sense we should consider the process $(X_{(s,t)})_{(s,t)\in\mathbb{R}_+^2}$ associated with the projective family $(\mu_{(s,t)})_{(s,t)\in\mathbb{R}_+^2}$ where
$$
\widehat{\mu}_{(s,t)}(\xi)=(2\pi)^{-\frac{n}{2}}e^{-s\psi_1(\xi)-t\psi_2(\xi)-d(s,t;\xi)}
$$
as an extension of the two L\'evy processes associated with $\big(T_s^{(1)}\big)_{s\geq 0}$ and $\big(T_t^{(2)}\big)_{t\geq 0}$ ``governed'' by the operator $a(s,t;D_x)$.
\end{enumerate}
This last example makes it also clear that we left the realm of two-(or multi-) parameter Markov processes, compare D. Khoshnevisan \cite{KHO} as standard reference.\\
\\
In some of the cases discussed above the interaction of probability theory and analysis were investigated, we mention again the work initiated by J. B. Walsh, E. B. Dynkin and others, and refer to Khoshnevisan. However we are missing any idea of a general theory, in particular when allowing operators with variable coefficients with respect to $x\in\mathbb{R}^n$. We may formulate the following general problem (2-parameter case only):
\begin{quote}
Given a pseudodifferential operator $a(s,t;x,D_x)$ with symbol $a:\mathbb{R}_+^2\times G\times \mathbb{R}^n$, $G\subset \mathbb{R}^n$. Find conditions under which the solutions to
\begin{align}\label{eq28}
\frac{\partial^2}{\partial s\partial t}v(s,t;x)=a(s,t;x,D_x)v(s,t;x)
\end{align}
give rise to a two parameter stochastic process and develop an appropriate stochastic analysis.
\end{quote}
We may take as starting point the problem
\begin{subequations}\label{eq29}
\begin{align}
\frac{\partial^2}{\partial s\partial t}v(s,t;x)=a(s,t;x,D_x)v(s,t;x), \label{eq29a}\\
\lim_{s\rightarrow 0}v(s,t;x)=T_t^{(2)}\varphi(x),\label{eq29b}\\
\lim_{t\rightarrow 0}v(s,t;x)=T_s^{(1)}\varphi(x),\label{eq29c}
\end{align}
\end{subequations}
where $\big(T_s^{(1)}\big)_{s\geq 0}$ and $\big(T_t^{(2)}\big)_{t\geq 0}$ are
Feller semigroups on $\mathbb{R}^n$. Note two special features of \eqref{eq29}. We do not specify data on the whole boundary or $(s,t)$-part of the boundary, but we are dealing with a distinguished boundary, namely $\{0\}\times\{0\}\times\mathbb{R}^n$ on which we prescribe data. Secondly, assume for a moment that
$a(s,t;x,\xi)$ is $x$ independent and that the semigroups $\big(T_s^{(1)}\big)_{s\geq 0}$ and $\big(T_t^{(2)}\big)_{t\geq 0}$ are spatially homogeneous, i.e. each is associated with a continuous negative definite function $\psi_j$. Then the Fourier transformed problem is
\begin{subequations}\label{eq30}
\begin{align}
\frac{\partial^2}{\partial s\partial t}\widehat{v}(s,t;\xi)&=a(s,t;\xi)\widehat{v}(s,t;\xi),\label{eq30a}\\
v(0,t;\xi)&=e^{-t\psi_2(\xi)}\widehat{\varphi}(\xi),\label{eq30b}\\[0.14cm]
v(s,0;\xi)&=e^{-s\psi_1(\xi)}\widehat{\varphi}(\xi).\label{eq30c}
\end{align}
\end{subequations}
Hence we face a hyperbolic problem with initial conditions on characteristic surfaces! In a forthcoming paper \cite{JACPOT} we start to investigate systematically problem \eqref{eq29}.

\subsection*{Acknowledgment}
The work of the second named author was supported by an RCUK-fellowship.

\vspace{1cm}
\noindent
Niels Jacob and Alexander Potrykus\\
Department of Mathematics\\
University of Wales Swansea\\
Singleton Park, Swansea\\
SA2 8PP\\
United Kingdom\\[0.3cm]
e-mail: \texttt{N.Jacob@swansea.ac.uk}\\
e-mail: \texttt{A.K.K.Potrykus@swansea.ac.uk}

\end{document}